\documentclass[12pt,bezier]{article}  

\usepackage{amsmath}
\usepackage{amssymb}
\usepackage{latexsym}
\usepackage{amssymb}
\usepackage{verbatim} 

\newif\ifdeveloping
\developingtrue

\ifdeveloping
\fi

\newcommand{\qed}{\rule{2.5mm}{3mm}\hfill\break\smallskip}
\newcommand{\ZZ}{\hbox{\boldmath{$Z$}}}
\newcommand{\Proof}{\noindent{\sc Proof. }}
\def\ZZ{{\hbox{\sf Z\kern-.43emZ}}}

\newtheorem{theorem}{Theorem}[section]
\newtheorem{corollary}[theorem]{Corollary}
\newtheorem{proposition}[theorem]{Proposition}
\newtheorem{lemma}[theorem]{Lemma}
\newtheorem{example}[theorem]{Example}
\newtheorem{definition}[theorem]{Definition}

\newtheorem{remark}[theorem]{Remark}
\newtheorem{result}[theorem]{Result}
\newtheorem{property}[theorem]{Property}

\title{Stability of $k$ mod $p$ multisets and small weight codewords of the code generated by the lines of PG$(2, q)$}
\author{Tam\'as Sz\H onyi 
and Zsuzsa Weiner \thanks{In the earlier phase of this research, both authors
were partially supported by OTKA Grant K 81310.  In the final phase, 
the first author was partially supported by the Slovenian-Hungarian OTKA Grant NN 114614.} }

\begin{document}
\maketitle

\begin{abstract}
In this paper, we prove a stability result on $k$ mod $p$ multisets of points in PG$(2, q)$, $q = p^h$. 
The particular case $k=0$ is used to describe small weight codewords of the code generated by the lines of PG$(2, q)$, as linear combination of few lines.
Earlier results proved this for codewords with weight less than $2.5q$, while our result is valid until $cq\sqrt q$.
It is sharp when $27 < q$ square and $h \geq 4$.  When $q$ is a prime, De Boeck and Vandendriessche (see \cite{{BoeckExample}}) constructed
a codeword of weight $3p-3$ that is not the linear combination of three lines. We characterise their example.
\end{abstract}

\section{Introduction}
In a previous paper (\cite{stabeven}), we proved a stability result on point sets of even type in PG$(2, q)$. 
A {\it set of even type} $S$ is a point set intersecting each line
in an even number of points.  It is easy to see that sets of even type can only exist when $q$ is even. 
A stability theorem says that when a structure is ``close'' to being extremal, then it can be obtained from 
an extremal one by changing it a little bit. More precisely, we proved that if the number of odd secants, $\delta$,  of a point set is less than
$(\lfloor \sqrt q \rfloor +1)(q+1-\lfloor \sqrt q \rfloor )$, then we can add and delete, altogether $\lceil {\delta\over q+1} \rceil$ points, so that we obtain a 
point set of even type. As a consequence, we described small weight codewords of $C_1(2, 2^h)$. 
 $C_1( 2, 2^h)$ is the binary code generated by the characteristic vectors of lines in PG$(2, 2^h)$. As the complement of an (almost) even set is an (almost) odd
 set, the same results hold for odd sets.
 
 The aim of this paper is to generalise the above results to odd $q$. A possible generalisation of sets of even type are sets intersecting every line in $k$ mod $p$ points,
 briefly {\it $k$ mod $p$ sets}.
 We expect the union of a $k_1$ mod $p$ and a $k_2$ mod $p$ set be  a $k_1 + k_2$ mod $p$ set. This is only true if we consider multisets, that is the points of the set have weights.
 We call a multiset  a {\it $k$ mod $p$ multiset} if it intersects every line in $k$ mod $p$ points counted with weights (multiplicities).
 Hence our aim is to generalise the stability results of sets of even type to $k$ mod $p$ multisets where $q = p^h$, $p$ prime.
 More precisely, the following theorems will be proved.
 
 \begin{theorem}
\label{linecover}
Let $\cal M$ be a multiset in ${\rm
PG}(2,q)$, $17< q$, $q=p^h$, where $p$ is prime.
Assume that the number of lines intersecting $\cal M$ in not
$k$ mod $p$ points is $\delta$, where $\delta <
{\sqrt {q \over 2} } (q+1)$. Then there
exists a set $S$ of points with size  $\lceil {\delta\over q+1}\rceil$, which blocks all the not 
$k$ mod $p$ lines.
\end{theorem}

\begin{theorem}
\label{kmodp} Let $\cal M$ be a multiset in ${\rm
PG}(2,q)$, $27< q$, $q=p^h$, where $p$ is prime and $h>1$ (that is $q$ not a prime).
Assume that the number of lines intersecting $\cal M$ in not
$k$ mod $p$ points is $\delta$, where 
\begin{itemize}
\item[(1)]{$\delta < (\lfloor \sqrt q \rfloor +1)(q+1-\lfloor \sqrt q \rfloor )$, when $2 < h$.}
\item[(2)]{$\delta < \frac{(p-1)(p-4)(p^2+1)}{2p-1}$, when $h = 2$.}
\end{itemize}
Then there
exists a multiset $\cal{M}'$ with the property that it intersects every
line in $k$ mod $p$ points and the number of different points in
 $({\cal M}\cup {\cal M}')\setminus ({\cal M}\cap {\cal M}')$ is
exactly $\lceil {\delta\over q+1}\rceil$.
\end{theorem}

Note that in $(1)$ $27< q$ and $2 < h$ actually means that $64 \leq q$ or $81 \leq q$ if q is odd.

\begin{remark}
{\em
Observe that the conclusion in Theorem \ref{kmodp} is much stronger than in Theorem \ref{linecover}, but Theorem 
\ref{kmodp} does not say anything when $h=1$ or $h=2$ and $\delta \geq \frac{(p-1)(p-4)(p^2+1)}{2p-1}$. Nevertheless, the conclusion in 
Theorem \ref{kmodp} does not apply in case $h=1$ as the result of De Boeck and Vandendriessche shows (see
Examples \ref{BelgaPelda} and  \ref{generalBelgaExample}.) }
\end{remark}


\begin{remark}
{\em
Note that a complete arc of size $q-\sqrt q+1$ has
$(\sqrt q+1)(q+1-\sqrt q)$ odd-secants, which shows that
Theorem \ref{kmodp} is
sharp, when $q$ is an even square. (Since the smallest sets of
even type are hyperovals.)
For the existence of such arcs, see \cite{BSz}, \cite{E}, \cite{FHT} and \cite{Ke}, \cite{W}.}
\end{remark}

Let $C_1(2,q)$ be the $p$-ary 
linear code generated by the characteristic vectors of the lines of ${\rm PG}(2,q)$
$q=p^h$, $p$ prime. Hence a codeword $c$ is a linear combination of lines, that is  $c = \sum_i \lambda_i l_i$.
The vectors $l_i$ correspond to a point in the dual plane of PG$(2, q)$. If we consider the point corresponding to $l_i$ with weight $\lambda_i$, then
$c$ corresponds to a multiset in the dual plane of PG$(2, q)$. A codeword $c$ with weight $w(c)$ (the number of non-zero coordinates) corresponds to a multiset
intersecting all but $w(c)$ lines in $0$ mod $p$ points. The coordinates of $c$ that are zero correspond to lines intersecting the multiset in the dual plane in $0$ mod $p$ points.
Hence we can translate our stability results on multisets (Theorem \ref{linecover} and \ref{kmodp}) to results on small weight codewords, 
see Theorem \ref{coverofsupp}, 4.3 and \ref{codewords}. Earlier results \cite{AK}, \cite{C} and \cite{FFSVW} described codewords of weight less than $2.5q$,
while our results are valid until $c  q\sqrt q$, when $q$ is not a prime. Here $c \approx \frac{1}{2}$, for $q=p^2$ and  $c \approx 1$ otherwise.
The prime case is a bit more difficult, because of some examples of weight $3p$ constructed by De Boeck and Vandendriessche, see \cite{{BoeckExample}} and Example \ref{BelgaPelda}. 
A slight generalisation that gives also codewords of weight $3p+1$ is given in Example \ref{generalBelgaExample}. 
In this case, we prove the following results.

\medskip 
\noindent {\bf Theorem 4.8} {\em Let $c$ be a codeword of $C(2,p)$, $p > 17$ prime.  If $2p+1 < w(c) \leq 3p+1$, then 
$c$ is either the linear combination of three lines or Example \ref{generalBelgaExample}. }

\medskip
\noindent {\bf Corollary \ref{intervals}} {\em For any integer $0< k + 1< {\sqrt{ q\over 2}}$, there is no codeword whose weight lies in the interval
 $(kq + 1, ((k+1)q -{3\over 2}k^2 - {5\over 2}k-1)$.}
 
 \medskip
 Note that, when $k = 3$, the above results give that codewords of weight less than $4p - 22$ can be obtained via Example \ref{generalBelgaExample} or they are linear combinations of
 three lines.

\section{The algebraic background}

In this section, we first recall the resultant-like condition formulated in \cite{stabeven} and \cite{W}.
Then we apply it for polynomials $g(X, Y)$ and $X^q-X$ and we state the combinatorial consequences
for the number of non-$k$ mod $p$ secants through a point.

\begin{result}
[\cite{W}, \cite{stabeven}]
\label{hn} 
Suppose that the nonzero polynomials $u(X,Y) = \sum_{i=0}^n u_i(Y)X^{n-i}$ and
$v(X,Y) = \sum_{i=0}^{n-m} v_i(Y)X^{n-m-i}$, $m>0$, satisfy {\em deg}$u_i(Y) \leq i$
and {\em deg}$v_i(Y) \leq i$ and $u_0 \not =0$.

Furthermore, assume that there exists a value $y$, so that the degree of the
greatest common divisor of $u(X,y)$ and $v(X,y)$ is $n-s$.
Denote by $n_h$, the number of values $y'$ for which
${\rm deg}({\rm gcd}(u(X,y'),v(X,y')))=n-(s-h)$.
Then $$\sum_{h=1}^{s}hn_h\leq s(s-m).\; \qed $$
\end{result}

\begin{remark}
{\em 
In our earlier paper \cite{stabeven}, unfortunately the index $h$ in Result \ref{hn}
ran until $s-1$ only, but the proof used $s$. However, the original Lemma 3.4 in \cite{W}
contains the right bound. Note that $h=s$ corresponds to $v = 0$.}
\end{remark}

We imagine the projective plane  PG$(2,q)$ as the closure of the affine plane AG$(2,q)$.
The affine points of PG$(2,q)$ will be represented by the pairs $(x,y)$.
 Affine lines with slope $m$ (with equation $Y = mX + b$)
pass through the point $(m)$, vertical lines (with equation $X = c$) pass through the infinite point $(\infty)$.

Let $\ell_{\infty}$ be the line at infinity intersecting the 
multiset ${\cal M}$ in $k$ mod $p$ points. Furthermore, let ${\cal
M}\setminus \ell_{\infty}=\{ (a_v,b_v)\}_v$ and ${\cal M}\cap
\ell_{\infty}=\{ (y_i) \}_i$, $(y_i) \not= (\infty)$. Consider the following polynomial:

\begin{equation}
g(X,Y)=\sum_{v=1}^{|{\cal M}\setminus
\ell_{\infty}|}(X+a_vY-b_v)^{q-1}+ \sum_{y_i\in {\cal M}\cap
\ell_{\infty}}(Y-y_i)^{q-1}-|{\cal M}| + k= \sum_{i=0}^{q-1}
r_i(Y)X^{q-1-i},
\end{equation}

\noindent Note that deg$r_i\leq i$.

\begin{lemma}
\label{index}
Through a point $(y)$ there pass $s$ non-$k$ mod $p$ affine secants
of $\cal M$ if and only if the degree of the greatest common divisor of
$g(X,y)$ and $X^q-X$ is $q-s$.
\end{lemma}

\Proof To prove this lemma, we only have to show that $x$ is a root
of $g(X,y)$ if and only if the line $Y=yX+x$ intersects $\cal M$ in
$k$ mod $p$ points.

Since $a^{q-1} = 1$, if $a \not = 0$ and $0^{q-1} = 0$, for the pair $(x, y)$
the number of zero terms in the first sum is exactly the number of affine points of $\cal M$
on the line $Y=yX+x$, the rest of the terms are $1$.  So assume that the ideal point $(y)$
of the line $\ell: Y=yX+x$ is in $\cal M$ with multiplicity $s$ ($0\leq s \leq p-1$).
Hence the first sum is $|{\cal M}| - k - (|\ell \cap {\cal M}| - s)$ (note that $|{\cal M} \cap l_{ \infty} | = k$).
The second sum is
$k-s$. Hence in total we get  $|{\cal M}| - k - (|\ell \cap {\cal M}| - s) + ( k- s) - |{\cal M}|  + k = k-|\ell \cap {\cal M}|$ 
and so the lemma follows. 
 \qed 


\begin{remark}
\label{odd:hn_h} {\em Assume that the line at infinity intersects ${\cal M}$
in $k$ mod $p$ points and suppose also that there is an ideal point,
different from $(\infty )$, with $s$ non-$k$ mod $p$ secants through it.
Let $n_h$ denote the number of ideal points
different from $(\infty)$, through which there pass $s-h$ non-$k$ mod $p$ secants
of the multiset ${\cal M}$. Then Lemma \ref{index} and Result
\ref{hn} imply that $\sum_{h=1}^{s} hn_h \leq  s(s-1)$. }
\end{remark}

\begin{lemma}
\label{maxindex}
Let ${\cal M}$ be a multiset in ${\rm PG}(2,q)$, $17<q$, 
so that the number of lines intersecting it in non-$k$ mod $p$
points is $\delta$, where  $\delta < (\lfloor \sqrt q \rfloor +1)(q+1-\lfloor \sqrt q \rfloor)$.
Then
\begin{enumerate}
\item[(1)]{the number $s$ of non-$k$ mod $p$ secants
through any point of  ${\cal M}$ satisfies $qs -s(s - 1) \leq \delta$,}
\item[(2)]{ the number of non-$k$ mod $ p$ secants through any point is at most
${\rm min}({\delta\over q+1}  + 2, \lfloor \sqrt q \rfloor + 1)$ or at least ${\rm max}(q + 1-({\delta\over q+1} + 2), q-\lfloor \sqrt q \rfloor)$.}
\end{enumerate}
\end{lemma}


\Proof
Pick a point $P$ with $s$ non-$k$ mod $p$ secants through it
and let $\ell_{\infty}$ be a $k$ mod $p$ secant of $\cal{M}$ through $P$. 
(If there was no such secant, then
the lemma follows immediately.)
By Remark \ref{odd:hn_h}, counting the number
of non-$k$ mod $p$ secants through the points of $\ell_{\infty} \setminus (\infty )$, we get $(1)$:
$$qs-s(s-1)\leq \delta .$$

The discriminant of the above inequality is $(q+1)\sqrt{1-\frac{4\delta}{(q+1)^2}}$. Let $x = \frac{4\delta}{(q+1)^2}$
and use $1 - {x\over 2} - {x^2\over 4} \leq \sqrt{1 - x}$ (which
is certainly true when $x \leq {4\over 5}$) to estimate the discriminant. In our case, $x  \leq{4\over 5}$, gives the condition $q > 17$. Hence $s \leq {\delta\over q+1} + {2\delta ^2\over (q+1)^3} (< {\delta\over q+1} + 2)$ or $s \geq q + 1-({\delta\over q+1} + {2\delta ^2\over (q+1)^3}) (> q-1 - {\delta\over q+1} )$.
On the other hand,  as the discriminant
is larger than $q+1-2(\lfloor \sqrt q \rfloor +2)$ (since $\delta < (\lfloor \sqrt q \rfloor +1)(q+1-\lfloor \sqrt q \rfloor)$), $s<\lfloor \sqrt q \rfloor +2$
or $s>q+1-(\lfloor \sqrt q \rfloor +2)$; whence $(2)$ of the lemma follows.
\qed

The next proposition is a generalisation of Lemma \ref{maxindex} and follows immediately.

\begin{proposition}
\label{generalindex}
Let ${\cal M}$ be a multiset in ${\rm PG}(2,q)$, $17<q$, 
so that the number of lines intersecting it in non-$k$ mod $p$
points is $\delta$, where  $\delta <  {3 \over 16} (q+1)^2$.
Then the number of non-$k$ mod $ p$ secants through any point is at most
${\delta\over q+1} + {2\delta ^2\over (q+1)^3}$ or at least $q + 1-({\delta\over q+1} + {2\delta ^2\over (q+1)^3})$. \qed
\end{proposition}

\section{ Proofs of Theorems \ref{linecover} and \ref{kmodp}}
\addtocounter{section}{-2}
 \begin{theorem}
\label{linecover}
Let $\cal M$ be a multiset in ${\rm
PG}(2,q)$, $17< q$, $q=p^h$, where $p$ is prime.
Assume that the number of lines intersecting $\cal M$ in not
$k$ mod $p$ points is $\delta$, where $\delta <
{\sqrt {q \over 2} } (q+1)$. Then there
exists a set $S$ of points with size  $\lceil {\delta\over q+1}\rceil$, which blocks all the not 
$k$ mod $p$ lines.
\end{theorem}

\begin{theorem}
\label{kmodp} Let $\cal M$ be a multiset in ${\rm
PG}(2,q)$, $27< q$, $q=p^h$, where $p$ is prime and $h>1$ (that is $q$ not a prime).
Assume that the number of lines intersecting $\cal M$ in not
$k$ mod $p$ points is $\delta$, where 
\begin{itemize}
\item[(1)]{$\delta < (\lfloor \sqrt q \rfloor +1)(q+1-\lfloor \sqrt q \rfloor )$, when $2 < h$.}
\item[(2)]{$\delta < \frac{(p-1)(p-4)(p^2+1)}{2p-1}$, when $h = 2$.}
\end{itemize}
Then there
exists a multiset $\cal{M}'$ with the property that it intersects every
line in $k$ mod $p$ points and the number of different points in
 $({\cal M}\cup {\cal M}')\setminus ({\cal M}\cap {\cal M}')$ is
exactly $\lceil {\delta\over q+1}\rceil$.
\end{theorem}

\noindent {\bf Proof of Theorem \ref{linecover}:}  
First we show that every line intersecting 
${\cal M}$ in non-$k$ mod $p$ points contains a point through which there are at least 
$q + 1-({\delta\over q+1} + {2\delta ^2\over (q+1)^3})$ lines intersecting ${\cal M}$ in non-$k$ mod $p$ points. We call such a
point a point with large index.
On the contrary, assume that $\ell$ is a line intersecting ${\cal M}$ in non-$k$ mod $p$ points but containing no point with large index. 
Then by Proposition \ref{generalindex},
through each point of $\ell$ there pass at most ${\delta \over q+1} + {2\delta ^2\over (q+1)^3}$ non-$k$ mod $p$ secants.
Hence $\delta \leq (q+1)({\delta\over q+1} + {2\delta ^2\over (q+1)^3} -1) + 1$. This implies that 
$\delta < {\sqrt{ q \over 2} } (q+1)$, a contradiction.

It is obviuos that to cover every line intersecting ${\cal M}$ in non-$k$ mod $p$ points we need at least 
$\lceil {\delta\over q+1} \rceil$ points.  We only need to show that there are less than ${\delta\over q+1}  + 1$ points with large index. Through every such point there are at least $q + 1-({\delta\over q+1} + {2\delta ^2\over (q+1)^3})$ non-$k$ mod $p$ secants, hence if there were at least ${\delta\over q+1}  + 1$ of them, then $$\delta \geq ({\delta\over q+1}  + 1)(q + 1-({\delta\over q+1} + {2\delta ^2\over (q+1)^3})) -{{\delta\over q+1} + 1 \choose 2}.$$
This is a contradiction since $\delta < {\sqrt{ q \over 2} } (q+1)$.\qed

\addtocounter{section}{2}

\begin{remark}
\label{blockingpointshavelargeindex}
{\em 
It follows from the beginning of the above proof that through each point of $S$ in Theorem \ref{linecover}, there
pass at least $q + 1-({\delta\over q+1} + {2\delta ^2\over (q+1)^3})$ lines intersecting ${\cal M}$ in non-$k$ mod $p$ points.}
\end{remark}

\begin{proposition}
\label{notthatmany}
Let ${\cal M}$ be a multiset in ${\rm PG}(2,q)$, $17<q$, having less
than $(\lfloor \sqrt q \rfloor +1)(q+1-\lfloor \sqrt q \rfloor)$
non-$k$ mod $p$ secants. Assume
that through each point there pass 
less than  $(q-\lfloor \sqrt q \rfloor )$
non-$k$ mod $p$ secants. Then the total number $\delta$ of lines intersecting
${\cal M}$ in non-$k$ mod $p$ points  is at
most $\lfloor \sqrt q \rfloor q-q+2\lfloor \sqrt q \rfloor + 1$.
\end{proposition}

\Proof Assume to the contrary that $\delta >\lfloor \sqrt q \rfloor
q-q+2\lfloor \sqrt q \rfloor + 1$. Pick a point $P$ and let
$\ell_{\infty}$ be a $k$ mod $p$ secant of ${\cal M}$ through $P$. Assume
that there are $s$ non-$k$ mod $p$ secants through $P$. If there is a point $Q$
on $\ell_{\infty}$ through which there pass at least $s$
non-$k$ mod $p$ secants, then choose the coordinate system so that $Q$ is
$(\infty )$. Then, by Remark \ref{odd:hn_h}, counting the number of
non-$k$ mod $p$ secants through $\ell$, we get a lower bound on $\delta$:

$$(q+1)s-s(s-1)\leq \delta .$$

Since $\delta < (\lfloor \sqrt q \rfloor +1)(q+1-\lfloor \sqrt q
\rfloor)$,  from the above inequality we get that $s<\lfloor \sqrt q
\rfloor +1$ (hence $s\leq \lfloor \sqrt q \rfloor$) or
$s>q+1-\lfloor \sqrt q \rfloor$, but by the assumption of the
proposition the latter case cannot occur.

Now we show that through each point there are at most
$\lfloor \sqrt q \rfloor$ non-$k$ mod $p$ secants.
The argument above and Lemma \ref{maxindex} show that on each $k$ mod $p$ secant
 there is at most one point through
which there pass $\lfloor \sqrt q \rfloor +1$ non-$k$ mod $p$ secants and through
the rest of the points there are at most $\lfloor \sqrt q \rfloor$ of them.
Assume that there is a point $R$
with $\lfloor \sqrt q \rfloor +1$ non-$k$ mod $p$ secants. Since
$\delta >\lfloor \sqrt q \rfloor +1$, we can find a non-$k$ mod $p$ 
secant $\ell$ not through $R$. From above, the number of non-$k$ mod $p$
secants through the intersection point of a $k$ mod $p$ secant on $R$ and $\ell$ is
at most $\lfloor \sqrt q \rfloor$. So counting the non-$k$ mod $p$ secants
through the points of $\ell$, we get at most
$(q-\lfloor \sqrt q \rfloor )(\lfloor \sqrt q \rfloor -1)+
(\lfloor \sqrt q \rfloor +1)\lfloor \sqrt q \rfloor +1$, which is a
contradiction. So there was no point with $\lfloor \sqrt q \rfloor +1$
non-$k$ mod $p$ secants through it.

This means that the non-$k$ mod $p$ secants form a dual $\lfloor \sqrt q \rfloor$-arc,
hence $\delta \leq (\lfloor \sqrt q \rfloor -1)(q+1)+1$, which is a
contradiction again; whence the proposition follows. \qed

\begin{property}
\label{most}
Let $\cal M$ be a multiset in ${\rm
PG}(2,q)$, $q=p^h$, where $p$ is prime. Assume that there are $\delta$
lines that intersect $\cal M$ in non-$k$ mod $p$ points. 
If through a point there are more than $q / 2$ lines intersecting 
 $\cal M$ in non-$k$ mod $p$ points, then there exists a value $r$ such that
 the intersection multiplicity of more than $2{\delta\over q + 1} + 5$ of these lines
 is $r$. 
 \end{property}
 
 Note that the property above does not hold necessary for all multisets having at most $\delta$ 
 non-$k$ mod $p$  secants. However, in Section 4, we are going to show that there are cases when the above property holds automatically.
 
 \medskip
 
 

\begin{theorem}
\label{general}
Let $\cal M$ be a multiset in ${\rm
PG}(2,q)$, $17< q$, $q=p^h$, where $p$ is prime.
Assume that the number of lines intersecting $\cal M$ in not
$k$ mod $p$ points is $\delta$, where $\delta <
(\lfloor \sqrt q \rfloor +1)(q+1-\lfloor \sqrt q \rfloor )$. 
Assume furthermore, that Property \ref{most} holds.
Then there
exists a multiset $\cal{M}'$ with the property that it intersects every
line in $k$ mod $p$ points and the number of different points in
 $({\cal M}\cup {\cal M}')\setminus ({\cal M}\cap {\cal M}')$ is
exactly $\lceil {\delta\over q+1}\rceil$.
\end{theorem}

Note that Theorem \ref{general} is also valid for $h=1$ and $h=2$ (not like Theorem \ref{kmodp}). Hence, for example, in case $h=1$ or $h=2$, if for a given set 
we know that Property \ref{most} holds, then Theorem \ref{general} yields a stronger result than Theorem \ref{kmodp}.

\medskip

\Proof
By Lemma
\ref{maxindex}, through each point there pass either at most ${\delta\over q+1} + 2$ 
or at least $q - 1- {\delta\over q+1} $ lines intersecting the multiset ${\cal M}$
in non-$k$ mod $p$ points. 
Let ${\cal P}$ be the set containing the points $P_i$ through which
there pass at least $q - 1- {\delta\over q+1} $ non-$k$ mod $p$ points. 
By Property \ref{most},  to each point $P_i$, there is a value $k_i$, so that more than
${\delta\over q+1} + 2$ lines through  $P_i$ intersect $\cal{M}$ in $k_i$ mod $p$ points. 
Add the point $P_1 \in {\cal P}$ to the multiset $\cal M$ with multiplicity $p-k_1$ and denote this new multiset
by ${\cal M}^{(1)}$. As there were only less than ${\delta\over q+1} + 2$ lines through $P_1$ which intersect
$\cal{M}$ in $k$ mod $p$ points and now by  Property \ref{most}, we ``repaired"   
more than ${\delta\over q+1} + 2$ lines, the total number of non-$k$ mod $p$ secants of ${\cal M}^{(1)}$ is less than $\delta$.
Hence again by Lemma
\ref{maxindex}, through each point there pass either at most ${\delta\over q+1} + 2$ 
or at least $q - 1- {\delta\over q+1} $ lines intersecting the multiset ${\cal M}^{(1)}$
in non-$k$ mod $p$ points. So, it follows that there are at most ${\delta\over q+1} + 2$
non-$k$ mod $p$ secants of ${\cal M}^{(1)}$ through $P_1$. It is also easy to see that the number of
non-$k$ mod $p$ secants of ${\cal M}^{(1)}$ is
at least $q + 1- 2( {\delta\over q+1} + 2 )$ less than that of $\cal M$.
Note that from the argument above and from Lemma
\ref{maxindex}, it also follows immediately that the  set containing the points through which
there pass at least $q - 1- {\delta\over q+1} $ non-$k$ mod $p$ secants of ${\cal M}^{(1)}$ is 
exactly ${\cal P} \setminus P_1$. We add the points of ${\cal P}$ one by one to ${\cal M}$
as above. At the $r$th step we want to add the point $P_r \in {\cal P}$ to ${\cal M}^{(r-1)}$. 
By Property \ref{most} and because of our algorithm, there are at least 
$2{\delta\over q + 1} + 5 - (r-1)$ lines through  $P_r$ intersecting ${\cal M}^{(r-1)}$ in $k_r$ mod $p$ points.
If $2{\delta\over q + 1} + 5 - (r-1)  > {\delta\over q+1} + 2$, then we can repeat the argument above and obtain the multiset
${\cal M}^{(r)}$. Note that at each step we ``repair" at least $q + 1- 2( {\delta\over q+1} + 2 )$ lines, hence there can be at most
${\delta \over  q + 1- 2( {\delta\over q+1} + 2 ) }$ steps in our algorithm, so our argument is valid at each step.

Let ${\cal M}'$ be the set which we obtain when ${\cal P}$ is empty and let $\delta '$
be the number of lines intersecting it in non-$k$ mod $p$ points. Proposition
\ref{notthatmany} applies and so  $\delta '\leq \lfloor \sqrt q \rfloor
q-q+2\lfloor \sqrt q \rfloor + 1$. 

{\it Our first aim is to show
that ${\cal M}'$ is a multiset intersecting each line in $k$ mod $p$ points.}

Let $P$ be an arbitrary point with $s$ secants intersecting ${\cal M}'$ in not
$k$ mod $p$ points, 
and let $\ell_{\infty}$ be a $k$ mod $p$ secant through $P$. Assume that
there is a point on $\ell_{\infty}$ with at least $s$ secants intersecting ${\cal M}'$
in non-$k$ mod $p$ points.
Then as in Proposition \ref{notthatmany}, counting the number of
non-$k$ mod $p$ secants through $\ell$, we get a
lower bound on $\delta '$:
$$(q+1)s-s(s-1)\leq \delta '.$$
This is a quadratic inequality for $s$, where the discriminant
is larger than  $(q+2-2{\delta '+q\over q+1})$.
Hence $s<{\delta '+q\over q+1}$ or $s>q+2-{\delta '+q\over q+1}$,
but by the construction of ${\cal M}'$, the latter case cannot occur.

Now we show that there is no point through which there pass
at least ${\delta '+q\over q+1}$ non-$k$ mod $p$ secants. On the contrary,
assume that $T$ is a point with ${\delta '+q\over q+1}\leq s$
non-$k$ mod $p$ secants. We choose our coordinate system
so that the ideal line is a $k$ mod $p$ secant through $T$ and
$T\not =(\infty )$.
Then from the argument above, through each ideal point,
there pass less than
$s(\geq {\delta '+q\over q+1})$ non-$k$ mod $p$ secants.
First we show that there exists an ideal point through which there
pass exactly $(s-1)$ non-$k$ mod $p$ secants.
Otherwise, by Remark \ref{odd:hn_h}, $2(q-1)\leq s(s-1)$;
but this is a contradiction since
$s\leq \lfloor \sqrt q \rfloor +1$ by Lemma  \ref{maxindex}. Let $(\infty )$ be a point
with $(s-1)$ non-$k$ mod $p$ secants. Then as before, we can give a lower
bound on the total number of non-$k$ mod $p$ secants of ${\cal M'}$:
$$(s-1)+qs-s(s-1)\leq \delta'$$
Bounding the discriminant (from below) by $(q+2-2{\delta '+q\over q+1})$,
it follows that $s<{\delta '+q\over q+1}$ or $s>q+2-{\delta '+q\over q+1}$. This is
a contradiction, since by assumption, the latter case cannot occur
and the first case contradicts our choice for $T$.

Hence through each point there
pass less than ${\delta '+q\over q+1}$ non-$k$ mod $p$ secants. Assume that $\ell$ is
a secant intersecting ${\cal M}'$ in non-$k$ mod $p$ points. Then summing up the non-$k$ mod $p$ secants
through the points of $\ell$ we get that
$\delta '<(q+1){\delta '-1\over q+1}+1$, which is a contradiction.
So ${\cal M}'$ is a multiset intersecting each line in $k$ mod $p$ points.

\medskip

To finish our proof we only have to show that the number of different points in
$({\cal M}\cup {\cal M}')\setminus ({\cal M}\cap {\cal M}')$ is $\lceil
{\delta\over q+1}\rceil $. As we saw in the beginning of this proof,
the number $\varepsilon $ of modified points
is smaller than $2\lfloor \sqrt q\rfloor$. On the one hand, if we
construct ${\cal M}$ from the set ${\cal M}'$ of $k$ mod $p$ type, then we
see that $\delta \geq \varepsilon (q+1-(\varepsilon -1))$. Solving the
quadratic inequality we get that $\varepsilon < \lfloor \sqrt q\rfloor +1$
or $\varepsilon >q+1-\lfloor \sqrt q\rfloor$, but from the argument above this latter
case cannot happen. On the other hand, $\delta \leq \varepsilon (q+1)$.
From this and the previous inequality (and from
$\varepsilon \leq \lfloor \sqrt q\rfloor $), we get that
${\delta \over q+1}\leq \varepsilon \leq {\delta \over q+1}+
{\lfloor \sqrt q\rfloor (\lfloor \sqrt q\rfloor -1)\over q+1}$. Hence the
theorem follows. \qed

\noindent{\bf Proof of Theorem \ref{kmodp}}
The previous proposition shows that 
to prove Theorem \ref{kmodp}, we only have to show that  Property \ref{most} holds.
By the pigeonhole principle, there is a value $r$, so that the intersection multiplicity of at least  $(q-1 - {\delta \over q + 1} )/(p-1)$ of the 
(non-$k$ mod $p$) lines with ${\cal M}$ is $r$. 
When $h > 2 $ and $q > 27$, then this is clearly greater than $ 2{\delta \over q+1} + 5$; hence Property \ref{most}  holds. In case
$h = 2$, assumption $(2)$ in the theorem ensures exactly that  $(p^2-1 - {\delta \over p^2 + 1} )/(p-1) > 2{\delta \over p^2+1} + 5 $ holds, so again
the property holds.
\qed

 \newpage
 \section{Codewords of PG$(2, q)$}
 
 \begin{definition}
Let $C_1(2,q)$ be the $p$-ary 
linear code generated by the incidence vectors of the lines of ${\rm PG}(2,q)$
$q=p^h$, $p$ prime. The weight $w(c)$ of a codeword $c$   $\in C_1(2,q)$ is the number of non-zero coordinates.
The set of coordinates, where $c$ is non-zero is denoted by supp$(c)$.
\end{definition}

The next theorem is a straightforward corollary of the dual of Theorem \ref{linecover}.

\begin{theorem}
\label{coverofsupp}
Let $c$ be a codeword of $C_1(2,q)$, with $17<q$, $q=p^h$, $p$ prime. If
$w(c) < \sqrt{{q \over 2 }} (q+1)$, then the points of
{\rm supp}$(c)$ can be covered by  $\lceil {w(c)\over q+1}\rceil $ lines.
\end{theorem}

\Proof 
By definition, $c$ is the linear combination of lines $l_i$ of PG$(2, q)$, that is $c = \sum_i \lambda_i l_i$.
For each point $P$, add the multiplicities $\lambda_i$ of the lines $l_i$ which pass through $P$.
By definition of the weight, there are exactly $w(c)$ points in 
PG$(2, q)$ through which this sum is not $0$ mod $p$. Hence the theorem follows from the dual of
Theorem \ref{linecover}.
\qed

Similary, from the dual of  Theorem \ref{kmodp}, we get the following theorem.

\begin{theorem}
Let $c$ be a codeword of $C_1(2,q)$, with $27<q$, $q=p^h$, $p$ prime. If
\begin{itemize}
\item{$w(c) <  (\lfloor \sqrt q \rfloor +1)(q+1-\lfloor \sqrt q \rfloor )$, $2 < h$, or }
\item{$w(c) <   \frac{(p-1)(p-4)(p^2+1)}{2p-1}$, when $h = 2$,}
\end{itemize}
then $c$ is a linear combination of exactly  $\lceil {w(c)\over q+1}\rceil $ different lines.
\end{theorem}

\Proof By definition, $c$ is a linear combination of lines $l_i$ of PG$(2, q)$, that is
$c = \sum_i \lambda_i l_i$. Let ${\cal C}$ be the multiset of lines where each line $l_i$ has multiplicity $\lambda_i$.
The dual of Theorem \ref{kmodp} yields that there are exactly 
 $\lceil {w(c)\over q+1}\rceil$ lines $m_j$ with some multiplicity $\mu_j$, such that
 if we add the lines $m_j$ with multiplicity $\mu_j$ to ${\cal C}$ then through any point of   PG$(2, q)$, 
 we see $0$ mod $p$ lines (counted with multiplicity). 
 
 In other words, we get that
$c + \sum_{j = 1}^{\lceil {w(c)\over q+1}\rceil }\mu_j m_j 
$  is the \underline 0 codeword.   
Hence $c = \sum_{j = 1}^{\lceil {w(c)\over q+1}\rceil }(p-\mu_j) m_j $.
\qed

Note that if we investigate proper point sets as codewords, then Property \ref{most} holds automatically. More precisely, let $B$ be a proper point set (each point has multiplicity $1$), which is
a codeword of $C_1(2, q)$. Hence $B$ corresponds to a codeword $c = \sum_i \lambda_i l_i$,  where $l_i$ are lines of PG$(2, q)$.
Again consider the dual of the multiset of lines where each line $l_i$ has multiplicity $\lambda_i$. Then, clearly, there are $w(c)$ lines intersecting this dual set in not $0$ mod $p$ point. Furthermore each of these lines has intersection multiplicity $1$ mod $p$ (as $B$ is a proper point set) and so Property \ref{most} holds; hence we can apply Theorem \ref{general}.

\begin{theorem}
\label{codewords}
Let $B$ be a proper point set in $PG(2, q)$, $17< q$. Suppose that $B$ is a codeword of the lines of $PG(2, q)$. Assume also that $|B| < (\lfloor \sqrt q \rfloor +1)(q+1-\lfloor \sqrt q \rfloor )$.
Then $B$ is the linear combination of at most $\lceil {|B| \over q+1} \rceil$ lines. \qed
\end{theorem} 

The following result summarises what was known about small weight codewords.

\begin{result}
\label{earliercodewordsresult}
Let $c$ be a non-zero codeword of $C_1(2,q)$, $q=p^h$, $p$ prime. Then
\begin{itemize}
\item[(1)]{{\em (Assmus, Key \cite{AK})} $w(c) \geq q + 1$. The weight of a codeword is  
$(q+1)$ if and only if the points corresponding to non-zero coordinates are the
$q+1$ points of a line.}
\item[(2)]{{\em (Chouinard \cite{C})} There are no codewords with weight in the closed interval $[q+2, 2q-1]$, for $h=1$.}
\item[(3)]{{\em (Fack, Fancsali, Storme, Van de Voorde, Winne \cite{FFSVW})} For $h=1$, the only codewords with weight at most 
$2p+(p-1) / 2$, are the linear combinations of at most two lines; so they have weight $p+1$, $2p$ or $2p+1$.  When $h > 1 $,
the authors exclude some values in the interval $[q + 2, 2q- 1]$. In particular, they exclude all weights in the interval $[ 3q/2, 2q -1 ]$, when $h \geq 4$.}
\end{itemize}  \qed
\end{result}

\begin{example}
\label{BelgaPelda} {\em (Maarten De Boeck, Peter Vandendriessche \cite{BoeckExample},  
Example 10.3.4)}
Let $c$ be a vector of the vector space $GF(p)^{p^2 + p + 1}$, $ p \not = 2$
a prime, whose positions correspond to the points of  {\em PG$(2, p)$}, such that
\[ c_P =\left\{
  \begin{array}{ll}
   a & if \  {P = (0, 1, a),}\\
   b & if \ {P = (1, 0, b),}\\
   c & if \ {P = (1, 1, c),}\\
    0 & \mbox{otherwise,}
  \end{array}
\right.
\] 

\noindent where $c_P$ is the value of $c$ at the position corresponding to the point 
$P$. Note that the points corresponding to positions with non-zero coordinates belong to the line 
$m : X_0 = 0$, the line $m' : X_1 = 0 $ or the line $m'' : X_0 = X_1$. These three lines are concurrent 
at the point $(0,0,1)$. Observe $w(c) = 3p-3$.
\end{example}

Next we generalise the example above. Note that, a collineation of the underlying plane PG$(2, q)$ induces a permutation on the
coordinates of $C_1(2,p)$, which maps codewords to codewords.

\begin{example}
\label{generalBelgaExample}
Let $c$ be the codeword in Example \ref{BelgaPelda}. Let $v_m$ be the incidence vector of the line $m$,
$v_{m'}$ the incidence vector of the line $m'$ and $v_{m''}$ of the line $m''$ in Example \ref{BelgaPelda}. 
Let $d:= \gamma c + \lambda v_m +  \lambda' v_m' + \lambda'' v_m''$.  Note that $w(d) \leq 3p+1$ as the points corresponding to
non-zero coordinates are on the three lines $m$, $m'$, $m''$.
Finally, let $\pi$ be a permutation on the 
coordinates induced by a projective transformation of the underlying plane {\em $PG(2,p)$}.
Our general example for codewords with weight at most $3p+1$ are the codewords $d$ with a permutation $\pi$ applied on its coordinate positions.
\end{example}


\begin{theorem}
\label{3p}
Let $c$ be a codeword of $C_1(2,p)$, $p > 17$ prime.  If $2p+1 < w(c) \leq 3p+1$, then 
$c$ is either the linear combination of three lines or given by Example \ref{generalBelgaExample}.
\end{theorem}

\Proof
By Theorem \ref{coverofsupp} (and since two lines can contain at most $2p+1$ points), supp$(c)$ can be covered by three lines $l_1, l_2, l_3$. 

 Assume that $c$ is in $C^{\perp}$ and the $l_is$ pass through the common point $P$. Note that as  $c$ is in $C^{\perp}$,
$P$  is  not in supp$(c)$. 
 First we show that either
 each multiplicity of the points in supp$(c)$ are different on each $l_i$, or the multiplicities of points of supp$(c)$ on a line $l_i$ are the same. 
 Let $S$ be the set of the points of  $l_1$ that have multiplicity $m$. Choose a point $Q$ from $l_2\setminus{ \{P \} }$, with multiplicity $m_Q$.
 As $c$ is in $C^{\perp}$, the multiplicities of the intersection points of any line with the $l_i$s should add
 up to $0$; hence the projection of $S$ from $Q$ onto $l_3$ is a set $S'$ of points with multiplicity $-( m_Q + m)$. Note that every point of $l_3$ outside $S'$
 must have multiplicity different from $-( m_Q+ m)$. Otherwise, projecting such a point back to $l_1$ from $Q$, the projection would have multiplicity $m$ (as $c$ is in the dual code); so it would be in $S$. Now pick a point $R$ of $l_3\setminus{ \{P \}}$  with multiplicity $n$ and choose a point $Q_R$, so that $Q_R$ projects to $R$ in $S$. From above, we see that there are exactly $|S|$ points on $l_3$ with multiplicity $n$ (which is the projection of $S$ from $Q_R$ onto $l_3$ ).
 This implies that $l_3\setminus{ \{P \}}$ is partitioned in sets of size $|S|$. As the number of points of $l_3\setminus{ \{P \}}$ is a prime, we get that $|S| = 1$ or $p$. 
 If the multiplicities of points of supp$(c)$ on a line $l_i$ are the same, then clearly $c$ is a linear combination of the lines $l_i$. 
 
We show that it is  Example \ref{generalBelgaExample}, when each multiplicity of the points in supp$(c)$ are different on each $l_i$. 
 As each point on $l_i$ has different multiplicity, let us choose our coordinate system, so that $P$ is the point $(0, 0, 1)$.
 The point of $l_1$ with multiplicity $0$ is the point $(0, 1, 0)$, the point of $l_3$ with multiplicity $0$ is the point $(1, 0, 0)$ and the point of 
 $l_2$ with multiplicity $-1$ is the point $(1, 1, 1)$.  Now we use the fact again that $c$ is in the dual code. 
 Hence from the line $[1, 0, 0]$ we get that the point $(0, 1, 1)$ has multiplicity $1$. 
 Examining line $[0, 1, 0]$ we get that the point $(1, 0, 1)$ has multiplicity $1$. Similarly if the point $(a, a, 1)$ has multiplicity $-m$,  we see that the points $(0, a, 1)$ and $(a, 0, 1)$ have 
 the same multiplicity, namely $m$. 
 Considering the line $<(0, 1, 1), (1, 0, 1)>$ we see that $(1/2, 1/2, 1)$ has multiplicity $-2$. So from above, the multiplicity of $(1 / 2, 0, 1)$ and  $(0, 1/2, 1)$ are $2$. 
 Now considering the line $<(0, 1 / 2, 1), (1 , 0, 1)>$ we see that  $(1/3, 1/3, 1)$ has multiplicity $-3$ and so $(1 / 3, 0, 1)$
and $(0, 1 / 3, 1)$ have multiplicity $3$. 
 Similarly, considering the line $<(0, 1 / n, 1), (1, 0, 1)>$ we see that  $(1/(n+1), 1/(n+1), 1)$ has multiplicity $-(n+1)$ and so $(1 / (n+1), 0, 1)$
and $(0, 1 / (n+1), 1)$ have multiplicity $(n+1)$; which shows that in this case $c$ is of Example \ref{generalBelgaExample}.
 
 Now assume that $c$ is in $C^{\perp}$, but the lines $l_i$ are not concurrent.  Assume that the intersection point $Q$ of $l_1\cap l_2$ has multiplicity $m$. Considering
 the lines through $Q$, we see that at least $(p-1)$ point on $l_3$ have multiplicity $-m$. Similarly, we see at least $(p-1)$ points on $l_2$
 and $(p-1)$ points on $l_3$ that have the same multiplicity. Hence taking the linear combination of $l_i$s with the right multiplicity, we get a codeword that only differs
 from $c$ in at most 3 positions (at the three intersection points of the lines $l_i$). There are no codewords with weight larger than 0 but at most $3$, which means that $c$ must be the linear combination of the $l_i$s. 

Now assume that $c$ is not in the dual code. As the dimension of the code is one larger than the dimension of the dual code (see \cite{H} and \cite{MM}), and $l_1$ is not in the dual code,
there exists a multiplicity $\lambda$, so that $c + \lambda l_1$ is in the dual code. It is clear that the weight of $c + \lambda l_1$ is  $\leq 3p$, and clearly supp$(c + \lambda l_1)$ can be covered by the three lines $l_1, l_2, l_3$. Now the result follows from the argument above when the weight of $c + \lambda l_1$ is greater than $2p$, and from Result \ref{earliercodewordsresult} otherwise. 
 \qed
 
 \begin{corollary}
 \label{intervals}
 For any integer $0< k + 1< {\sqrt{ q\over 2}}$, there is no codeword whose weight lies in the interval
 $(kq + 1, ((k+1)q -{3\over 2}k^2 - {5\over 2}k-1)$, for $q > 17$.
 \end{corollary}
 
 \Proof 
 Suppose to the contrary that $c$ is a codeword whose weight lies in the interval
 $(kq + 1, ((k+1)q -{3\over 2}k^2 - {5\over 2}k-1)$. Then by Theorem \ref{coverofsupp}, supp$(c)$ can be covered by the set $k+1$
 lines $l_i$. It follows from Remark \ref{blockingpointshavelargeindex}, that the number of points of supp$(c)$ on a line $l_i$ is at least $q - k - 1$.
 Hence $w(c)$ is at least $(k+1)(q - k - 1)- \binom{k+1}{2}$.
 
 
 \begin{corollary}
 \label{charofcodewords}
 Let $c$ be a codeword of $C(2,p)$, $p > 17$ prime.  If  $w(c) \leq 4p - 22$, then 
$c$ is either the linear combination of at most three lines or Example \ref{generalBelgaExample}.
 \end{corollary}
 
 \Proof It follows from Corollary \ref{intervals}, Theorem \ref{3p} and Result \ref{earliercodewordsresult}. \qed
 
 \bigskip
 \noindent
{\bf Acknowledgment.}
 The results on small weight codewords were inspired by conversation with Andr\'as G\'acs. 
 We gratefully dedicate this paper to his memory.

\noindent Authors address:\\
\noindent
Tam\'as Sz\H onyi\\
Department of Computer Science, E\"otv\"os Lor\'and University,\\
H-1117 Budapest, P\'azm\'any P\'eter s\'et\'any 1/C, HUNGARY\\
\noindent {\tt e-mail: szonyi@cs.elte.hu}\\

\noindent
Tam\'as Sz\H onyi, Zsuzsa Weiner\\
MTA-ELTE Geometric and Algebraic Combinatorics Research Group,\\
 H-1117 Budapest, P\'azm\'any P\'eter s\'et\'any 1/C, HUNGARY\\
 \noindent {\tt e-mail: zsuzsa.weiner@gmail.com}\\
 
 \noindent
 Zsuzsa Weiner\\
 Prezi.com\\
 H-1065 Budapest, Nagymez\H o utca 54-56, HUNGARY\\
 \end{document}